\documentclass{amsart}
\usepackage{amsmath}
\usepackage{amssymb}

\def\xx{{\bf x}}
\def\ww{{\bf w}}

\def\uu{{\bf u}}
\def\vv{{\bf v}}
\def\+{{\oplus}}

\newcommand{\BBZ}{\mathbb{Z}}

\theoremstyle{plain}
\newtheorem{theorem}{Theorem}
\newtheorem{corollary}{Corollary}
\newtheorem{lemma}{Lemma}

\newtheorem{remark}{Remark}

\theoremstyle{definition}
\newtheorem{definition}{Definition}

\theoremstyle{remark}

\begin{document}
\title[]{ Balanced Symmetric Functions over $GF(p)$}
\author{Thomas W. Cusick$^1$, Yuan Li$^2$,  Pantelimon St\u anic\u a$^{3*}$}
\address{$^1$SUNY, Department of Mathematics\\
         244 Mathematics Building\\
         Buffalo, NY 14260}
\email{email: cusick@buffalo.edu}
\address{$^2$Department of Mathematical Sciences\\
Alcorn State University\\
Alcorn State, MS 39096}
\email{email: yuanli7983@gmail.com}
\address{$^3$Applied Mathematics Department\\
Naval Postgraduate School\\
Monterey, CA 93943}
\email{email: pstanica@nps.edu}
\thanks{$^*$ Research supported by the Naval Postgraduate School RIP funding.}

\date{\today}

\begin{abstract}
  Under mild conditions on $n,p$, we give a lower bound on the number of $n$-variable balanced symmetric
polynomials  over finite fields $GF(p)$, where $p$ is a prime number.
The existence of nonlinear balanced symmetric polynomials is an
immediate corollary of this bound.
Furthermore, we conjecture that $X(2^t,2^{t+1}l-1)$ are the only nonlinear balanced
elementary symmetric polynomials over $GF(2)$, where
$X(d,n)=\sum_{i_1<i_2<\cdots<i_d}x_{i_1}x_{i_2}\cdots x_{i_d}$, and
we prove various results in support of this conjecture.
\end{abstract}
\maketitle

\section{Introduction}
\label{S:intro}
Since symmetry guarantees that all of the input bits have equal status in a very strong sense,
symmetric Boolean functions display some interesting properties. A lot
of research about symmetry in characteristic 2 has been previously done in
\cite{ann, carlet,Wu1,ssac,Wu2,Stinson,Gath,Sub,Mitchell,Sarkar,PS,Yuan4,Yang}. On the other hand,
it is natural to extend various
cryptographic ideas from $GF(2)$ to other finite fields of characteristic $>2$,
$GF(p)$ or $GF(p^n)$, $p$ being a prime number. For example, \cite{Mulan} and
\cite{Xiao} studied the correlation immune and resilient functions on $GF(p)$. Also, \cite{Feng} and
\cite{Kumar} investigated the generalized bent functions on $GF(p^n)$. In \cite{Yuan}, Li and Cusick first
introduced the strict avalanche criterion over $GF(p)$. In \cite{Yuan4}, they generalized most
 results of \cite{Wu2} and determined all the linear structures of symmetric functions over $GF(p)$.

 Balancedness is a desirable requirement of functions which will be used in cryptography. In
 this paper,  by an enumerating method, we give a lower bound for the number
 of balanced symmetric polynomials
 over $GF(p)$, and as an immediate consequence, we show the existence of nonlinear
 balanced symmetric polynomials.
 We did not find (even conjecturally) any simple characterization of the
algebraic normal form of nonlinear balanced symmetric polynomials even
for $p=2$.  However, we do make substantial progress in the binary case if
the polynomial is elementary symmetric (Section 5 below).
 We prove some results toward the conjecture
 that the polynomials $X(2^t,2^{t+1}\ell-1)$ are the {\em only} nonlinear balanced elementary symmetric polynomials,
 where $X(d,n)=\sum_{i_1<i_2<\cdots<i_d}x_{i_1}x_{i_2}\cdots x_{i_d}$.

\section{Preliminaries} \label{Pre}
In this paper, $p$ is a  prime number.
If $f$: $GF(p)^n\longrightarrow GF(p)$, then $f$ can be uniquely expressed in the following
form, called the {\em algebraic normal form} (ANF):
\[
f(x_1,x_2,\ldots,x_n)=\sum_{k_1,k_2,\ldots,k_n=0}^{p-1}
a_{k_1k_2\ldots k_n}{x_1}^{k_1}{x_2}^{k_2}\cdots {x_n}^{k_n},
\]
where each coefficient $a_{k_1k_2\ldots k_n}$ is a constant in $GF(p)$.

The function $f(x)$ is called an {\em affine function} if $f(x)=a_1x_1+\cdots+a_nx_n+a_0$.
If $a_0=0$, $f(x)$ is also called a {\em linear function}. We will denote by $F_n$ the set of all
 functions of $n$ variables and by $L_n$ the set of affine ones. We will
call a function {\em nonlinear} if it is not in $L_n$.

If $f(x)\in F_n$, then $f(x)$ is a {\em symmetric function} if for any permutation $\pi$
on $\{1, 2,\ldots,n\}$, we have $f(x_{\pi(1)}, x_{\pi(2)}, \ldots,x_{\pi(n)})$=$f(x_1, x_2,\ldots,x_n)$.
The set of permutations on $\{1, 2,\ldots,n\}$ will be denoted by $S_n$.

We define the following equivalence relation on $GF(p)^n$:
for any $x=(x_1,  \ldots,x_n)$, $y=(y_1,  \ldots,y_n)$ in $GF(p)^n$, we say $x$ and $y$
are equivalent, and write $x\sim y$, if there exists a permutation $\pi\in S_n$
such that $(y_1, y_2,\ldots, y_n)$=$(x_{\pi(1)}, x_{\pi(2)},\ldots,x_{\pi(n)})$    (by
abuse of notation we write $y=\pi(x)$). Let $\widetilde{x}=
\{y\,|\, \exists\, \pi\in S_n, \pi(x)=y\}$. Let
$\overline{x}=(\overline{x_1}, \overline{x_2},\ldots, \overline{x_n})$ be the
representative of $\widetilde{x}$, where $0\leq \overline{x_1}\leq
\overline{x_2}\leq\cdots\leq \overline{x_n}\leq p-1$.
Obviously, we have $\widetilde{x}=\widetilde{y}$ $\iff$ $\overline{x}=\overline{y}$.

\section{Enumeration Results}\label{Enumeration}

\begin{definition}
\label{balance-defn}
$f$: $GF(p)^n\longrightarrow GF(p)$ is {\em balanced} if
the probability $prob(f=k)=\frac{1}{p}$ for any $k=0,1,\ldots,p-1$.
\end{definition}
As an immediate consequence, $f$ is balanced if and only if $\#\{x\in GF(p)^n|f(x)=k\}=p^{n-1}$.


Using the equivalence relation of the previous section, we get that
$f$: $GF(p)^n\longrightarrow GF(p)$ is  symmetric if $f(x)=f(y)$ whenever $\widetilde{x}=\widetilde{y}$.
Let $C(n, k)=\frac{n!}{k!(n-k)!}$ if $0\leq k\leq n$ and $0$ otherwise
be the usual binomial coefficients.  Then we have
\begin{lemma}
\label{symmetry}
The number of $n$-variable symmetric polynomials over $GF(p)$ is
\[
p^{C(p+n-1, n)}.
\]
\end{lemma}
\begin{proof}
The number of different vector classes $\widetilde{x}$ is the number of solutions
of the linear equation $i_0+i_1+\cdots+i_{p-1}=n$, where
$i_k$ is the number of times  $k$ appears in $\overline{x}$. We know that the number of solutions
to the previous linear diophantine equation is the same as the number
of $n$-combinations of a set with $p$ elements, that is
$C(p+n-1, n)$ (see \cite[p. 69]{book}). Since a symmetric function $f(x)$ has the same value for any element of
$\widetilde{x}$, the lemma is proved.
\end{proof}
\begin{lemma}\label{3.2}
We have
$\displaystyle \prod_{k=0}^{p-1}C((k+1)a,a)=\frac{(pa)!}{(a!)^p}$.
\end{lemma}
\begin{proof}
It is a  straightforward computation
\[\prod_{k=0}^{p-1}C((k+1)a,a)=\frac{a!}{a!}\frac{(2a)!}{a!a!}
\frac{(3a)!}{a!(2a)!}\cdots\frac{(pa)!}{a!((p-1)a)!}= \frac{(pa)!}{(a!)^p}.
\]
\end{proof}
\begin{lemma}
\label{balance}
The number of $n$-variable balanced polynomials over $GF(p)$ is
\begin{displaymath}
\frac{(p^n)!}{(p^{n-1}!)^p}.
\end{displaymath}
\end{lemma}
\begin{proof}
The number we are looking for is
\[
C(p^n, p^{n-1})C(p^n-p^{n-1}, p^{n-1})\cdots
C(p^n-(p-1)p^{n-1}, p^{n-1})=\frac{(p^n)!}{(p^{n-1}!)^p},
\]
using Lemma \ref{3.2}, and the claim is proved.
\end{proof}

Let $\overline{x}=(\underbrace{0,\ldots, 0}_{i_0},
\underbrace{1,\ldots, 1}_{i_1},\ldots,\underbrace{p-1,\ldots, p-1}_{i_{p-1}}),$
where $i_0+i_1+\cdots+i_{p-1}=n$, $0\leq i_j\leq n$, $j=0, 1,\ldots,p-1$. The cardinality of the set
$\widetilde{x}$ is the value of the multinomial coefficient $C(n,i_0,i_1,\ldots,i_{p-2})=
\frac{n!}{i_0!i_1!\cdots i_{p-1}!}$.
We have the following widely known multinomial expansion lemma.
\begin{lemma} \cite[p. 123]{book}
\label{3.3}
We have the following formula
\begin{displaymath}
(t_0+t_1+\cdots+t_{p-1})^n=\sum_{\substack{i_0+i_1+\cdots+i_{p-1}=n}}C(n, i_0, i_1,\ldots,
i_{p-2})t_0^{i_0}t_1^{i_1}\cdots t_{p-1}^{i_{p-1}}.
\end{displaymath}
\end{lemma}
By specializing  $t_0=t_1=\cdots =t_{p-1}=1$, we get the following corollary.
\begin{corollary}
\label{3.4}
The $n$-th power of $p$ satisfies
\begin{displaymath}
p^n=\sum_{\substack{i_0+i_1+\cdots+i_{p-1}=n}}C(n, i_0, i_1,\ldots, i_{p-2}).
\end{displaymath}
\end{corollary}

From the proof of Lemma \ref{symmetry}, we know that
the number of terms in the sum in Corollary \ref{3.4} is $C(p+n-1,n)$. It is clear now, that
to get balanced symmetric polynomials amounts to partitioning the set of
$C(p+n-1,n)$ many multinomial coefficients $C(n,i_0,i_1,\ldots,i_{p-2})$ into $p$ groups,
the sum of each group being equal to $p^{n-1}$.

For a fixed solution $\{i_0,i_1,\ldots,i_{p-1}\}$ of $i_0+i_1+\cdots +i_{p-1}=n$, there are
 $\frac{p!}{m_0!m_1!\cdots m_n!}$ many ways to order it, where
 $i_j\in \{0,1,\ldots ,n\}$, and $m_l$ is the number of times that $l$ appears in
 $\{i_0,\ldots,i_{p-1}\}$, $0\le l\le n$. Hence,
 \begin{equation}
 \label{eqaaa}
 m_0+m_1+\cdots +m_n=p,\
\text{and}\  0 m_0+1m_1+\cdots +n m_n=n.
\end{equation}
Let us consider the following map:
\[
F: \{\{i_0,i_1,\ldots ,i_{p-1}\}|\sum_{j=0}^{p-1}i_j=n\}\to
 \{(m_0,m_1,\ldots ,m_n)|\sum_{l=0}^{n}m_l=p,\sum_{l=0}^{n}lm_l=n\}
 \]
defined by
\begin{displaymath}
F(\{i_0,i_1,\ldots,i_{p-1}\})=(m_0,m_1,\ldots,m_n),
\end{displaymath}
where $m_l$ is as above.
It is not hard to check that $F$ is a bijection.

Now, we will partition the set
of  multinomial coefficients $C(n, i_0,\ldots, i_{p-2})$ using the following
equivalence relation:
$C(n, i_0,\ldots , i_{p-2})$ and $C(n, j_0,\ldots, j_{p-2})$ belong to the
same class if and only if $j_0,\ldots,j_{p-1}$ is a permutation of
$i_0,\ldots,i_{p-1}$. Of course, any element in the same class has the same value.
So, we can think of $F$ as a map that assigns to each class the value
$\frac{p!}{m_0!m_1!\cdots m_n!}$.
\begin{lemma}\label{3.5}
Let $n,p$ be positive integers, with
$p$ a prime number. If $m_i<p$ for some $i$ (and so for all $i$), or if $\gcd(n,p)=1$, then
 $p$ divides $\frac{p!}{m_0!m_1!\cdots m_n!}$.
\end{lemma}
\begin{proof}
Assume $m_i<p$. By a known extension of Kummer's result that belongs to  Dickson
(see \cite[Theorem D, p. 3860]{Holte})
the power of $p$ that divides the multinomial coefficient equals the number of carries
when we add $m_0+m_1+\cdots+m_n$ in base $p$, but the mentioned sum is equal to $p$, therefore
the number of carries is~1. (One can also prove the same assertion without using Dickson's result.)

Now, assume $\gcd(n,p)=1$.
If $m_i<p$,  the first part of the proof proves the claim. Assume $m_i\geq p$. Since
$m_0+m_1+\cdots +m_n=p$, we can find $j$ such that $m_j=p$ and
$m_0=\cdots =m_{j-1}=m_{j+1}=\cdots m_n=0$.
From the definition of the $m_i$'s we obtain that $jp=n$, which is a contradiction.
\end{proof}

\begin{remark}
The two conditions $m_i<p$, and $\gcd(n,p)=1$ are not equivalent (although,
it is true that $\gcd(n,p)=1$ implies $m_i<p$). For instance, by taking
$m_0=3,m_1=2,m_2=1,m_3=1,m_4=m_5=m_6=m_7=0$, we get $m_0+m_1+\cdots +
m_7=p=7=n=0m_0+1m_1+\cdots+7m_7$,
so $p=n$ in this case.
\end{remark}

Since the cardinality of each multinomial coefficient class is a multiple of $p$, we can divide
each class into $p$ groups with an equal number of coefficients, hence, equal sum. Doing the same
for each class, we finally partition all of the $C(p+n-1,n)$ coefficients into $p$ groups with equal sum.

For a given $(m_0,m_1,\ldots ,m_n)$, $m_0+m_1+\cdots +m_n=p$, $0m_0+1m_1+\cdots +nm_n=n$,
the partition number is
{\tiny
\begin{eqnarray*}
&& C\left(\frac{p!}{m_0!m_1!\cdots m_n!},\frac{(p-1)!}{m_0!m_1!\cdots m_n!}\right)
 C\left(\frac{p!}{m_0!m_1!\cdots m_n!}-\frac{(p-1)!}{m_0!m_1!\cdots m_n!},
\frac{(p-1)!}{m_0!m_1!\cdots m_n!}\right)\cdots\\
&& C\left(\frac{p!}{m_0!m_1!\cdots m_n!}-\frac{k(p-1)!}{m_0!m_1!\cdots m_n!},
\frac{(p-1)!}{m_0!m_1!\cdots m_n!}\right) \cdots C\left(\frac{(p-1)!}{m_0!m_1!\cdots m_n!},
\frac{(p-1)!}{m_0!m_1!\cdots m_n!}\right).
\end{eqnarray*}
}
By Lemma \ref{3.2}, this product can be written as
\begin{displaymath}
\frac{(\frac{p!}{m_0!\cdots m_n!})!}{((\frac{(p-1)!}{m_0!\cdots m_n!})!)^p}.
\end{displaymath}
In conclusion, we get our main result of this section.
\begin{theorem}\label{1}
Let $N$ be the number of $n$-variable balanced symmetric functions over $GF(p)$.
If $m_i<p$, for all $i$ (or $\gcd(n,p)=1$), then
\begin{displaymath}
N\geq \prod_{\substack {\sum_{j=0}^{n}m_j=p\\\sum_{j=0}^{n}jm_j=n}}
\frac{(\frac{p!}{m_0!\cdots m_n!})!}{((\frac{(p-1)!}{m_0!\cdots m_n!})!)^p}.
\end{displaymath}
\end{theorem}

Next, since the linear balanced symmetric polynomials over $GF(p)$ have the form
$a(x_1+\cdots +x_n)+b$,
where $a\in GF(p)^*$ and $b\in GF(p)$, we get that the number of such functions is $p(p-1)$.
Since $\frac{(pa)!}{(a!)^p}=\frac{a!}{a!}\frac{(2a)!}{a!a!}\frac{(3a)!}{a!(2a)!}\cdots
\frac{(pa)!}{a!((p-1)a)!}>12\cdots p=p!\geq p(p-1)$, we have the next corollary.
\begin{corollary}\label{col}
If $n$ is not divisible by $p$, there exists a nonlinear $n$-variable balanced symmetric polynomial over $GF(p)$.
\end{corollary}


\section{The balancedness of elementary symmetric polynomials over $GF(2)$}

In this section we consider the binary case, that is, $p=2$. 
Here, we shall try to find all nonlinear balanced elementary symmetric  polynomials.
Throughout, $\xx=(x_1,\ldots,x_n)$ and $\oplus$ is the addition modulo 2.
\begin{definition}\label{de1}
For integers $n$ and $d$, $1\le d\le n$ we define the elementary symmetric polynomial by
\begin{equation}
\label{eq:defn-sigma}
X(d,n)=\sum_{i_1<i_2<\cdots<i_d}x_{i_1}x_{i_2}\cdots x_{i_d}.
\end{equation}
\end{definition}
By abuse of notation, we let $X(d,n)(j)$ be the value of $X(d,n)$
when $wt(\xx)=j$. Since $X(d,n)(j)\equiv C(j,d)\pmod{2}$, we get
\[
X(d,n)(j)=\frac{1-(-1)^{C(j,d)}}{2}.
\]
Because there are $C(n,j)$ many vectors with weight $j$, we have the following theorems.
\begin{theorem}
\label{th2}
The elementary symmetric polynomial $X(d,n)$ is balanced if and only if
\[
\sum_{0\le j\le n}C(n,j)(-1)^{C(j,d)}=0.
\]
\end{theorem}
\begin{theorem}\label{th3}
If $X(d,n)$ is balanced, then $d\le \lceil{n/2}\rceil$.
\end{theorem}
\begin{proof}
 If $n$ is even and $d\ge \frac{n}{2}+1$, then
\[
  \sum_{C(j,d)\equiv 0\pmod 2}C(n,j)> C(n,0)+C(n,1)+\cdots+C(n,n/2)>2^{n-1}.
\]

  If $n$ is odd and $k\ge \frac{n+1}{2}+1$, then
\[
  \sum_{C(j,d)\equiv 0\pmod 2}C(n,j)> C(n,0)+C(n,1)+\cdots+C(n,(n+1)/2)>2^{n-1}.
  \]

  In both cases, we have
\begin{eqnarray*}
&& \sum_{0\le j\le n}C(n,j)(-1)^{C(j,d)}\\
&&=\sum_{C(j,d)\equiv 0\pmod 2}C(n,j)-
\sum_{C(j,d)\equiv 1\pmod 2}C(n,j)\\
&&=
\sum_{C(j,d)\equiv 0\pmod 2}C(n,j)-\left(2^n-\sum_{C(j,d)\equiv 0\pmod 2}C(n,j)\right)\\
&&=
2\left(\sum_{C(j,d)\equiv 0\pmod 2}C(n,j)-2^{n-1}\right)>0,
\end{eqnarray*}
 contradicting Theorem \ref{th2}.
\end{proof}

Therefore, we see from Theorem \ref{th2} that the existence of balanced elementary
symmetric polynomials
is related to the problem of bisecting binomial coefficients (defined below).
In \cite{ssac}, two of us found some computational results about such bisections, which results
we shall describe below.
 (We mention here that the authors of \cite{Sarkar} found the number
 of solutions but without the explicit solutions.)
It was suspected that the existence of nontrivial binomial coefficient bisections (as in
\cite{ssac}) may cause
difficulties in the study of the existence of balanced symmetric polynomials, but we conjecture that
     this is not true for the elementary symmetric case.

 We begin with
\begin{definition}\label{de4}\cite{ssac}
If $\sum_{i=0}^{n} \delta_i C(n,i)=0,\ \delta_i\in\{-1,1\},\ i=0,1,\ldots,n$,
we call $(\delta_0,\ldots ,\delta_{n})$ a solution of the equation
\begin{equation}
\label{eq3}
\sum_{i=0}^{n}x_iC(n,i)=0,\quad x_i\in\{-1,1\}.
\end{equation}
\end{definition}

In fact, whenever we get a solution of (\ref{eq3}), we get a bisection of binomial
coefficients, that is, we find $A$, $B$ such that $A\cup B=\{0,1,\ldots ,n\}$,
$A\cap B=\emptyset$,
$\sum_{i\in A}C(n,i)=\sum_{i\in B}C(n,i)=2^{n-1}$.

Obviously, if $n$ is even, then $\pm(1,-1,1,-1,\ldots,1)$ are two solutions of (\ref{eq3}).
If $n$ is odd, then $(\delta_0,\ldots,\delta_{\frac{n-1}{2}},
-\delta_{\frac{n-1}{2}-1},\ldots,-\delta_0)$
are $2^\frac{n+1}{2}$ solutions of (\ref{eq3}). We call these trivial solutions.

Mitchell \cite{Mitchell} mentioned the nontrivial solutions for $n=8,13$.
In \cite{ssac}, with a C++ program,
we found all solutions of (\ref{eq3}) when $n\leq 28$. Nontrivial
solutions exist if and only if $n=8,13,14,20,24,26$.
So, here we ask the question of determining necessary and sufficient conditions
on the parameter
$n$ such that there exist nonlinear balanced symmetric polynomials on $GF(2)^n$.

First, we recall a known result that enables one to find residues of binomial
coefficients modulo a prime $p$.
\begin{lemma}[Lucas' Theorem]
\label{le7}
Let $n=a_mp^m+a_{m-1}p^{m-1}+\cdots +a_1p+a_0$ with $0\le a_i\le p-1$ and
$k=b_mp^m+b_{m-1}p^{m-1}+\cdots +b_1p+b_0$ with $0\le b_i\le p-1$, then
$C(n,k)\equiv C(a_m,b_m)\cdots C(a_1,b_1) \pmod{p}$
\end{lemma}
The next lemma can be derived from \cite{ann}. However, here we give a direct proof.
\begin{lemma}\label{le8}
For any integer $d\ge 2$, the sequence $\{(-1)^{C(j,d)}\}_{j=0}^{\infty}$ is
periodic of least period
$2^{[\log_2d]+1}$.
\end{lemma}
\begin{proof}
First, recall that $d$ has at most $[\log_2d]+1$ bits. For $0\le i\le 2^{[\log_2d]+1}-1$,
according to Lemma \ref{le7}, we have
$C(i+2^{[\log_2d]+1},d)\equiv C(1,0)C(i,d)\equiv C(i,d)\pmod{2}$, so
the least period is a divisor
of $2^{[\log_2d]+1}$. On the other hand,
$1=C(d,d)$ and $C(d+2^{[\log_2d]},d)\equiv C(1,0)C(0,1)\cdots \equiv 0 \pmod{2}$,
which implies that
$2^{[\log_2 d]}$ cannot be a period. The lemma is proved.
\end{proof}

With the help of Lemma \ref{le8}, we get the following computational results.
The list could easily be extended.
\begin{lemma}\label{le9}
We have\\
\indent
$\{\frac{1-(-1)^{C(j,2)}}{2}\}_{j=0}^{\infty}=\overline{0011}$

$\{\frac{1-(-1)^{C(j,3)}}{2}\}_{j=0}^{\infty}=\overline{0001}$

$\{\frac{1-(-1)^{C(j,4)}}{2}\}_{j=0}^{\infty}=\overline{00001111}$

$\{\frac{1-(-1)^{C(j,5)}}{2}\}_{j=0}^{\infty}=\overline{00000101}$

$\{\frac{1-(-1)^{C(j,6)}}{2}\}_{j=0}^{\infty}=\overline{00000011}$

$\{\frac{1-(-1)^{C(j,7)}}{2}\}_{j=0}^{\infty}=\overline{00000001}$

$\{\frac{1-(-1)^{C(j,8)}}{2}\}_{j=0}^{\infty}=\overline{0000000011111111}$

$\{\frac{1-(-1)^{C(j,9)}}{2}\}_{j=0}^{\infty}=\overline{0000000001010101}$

$\{\frac{1-(-1)^{C(j,10)}}{2}\}_{j=0}^{\infty}=\overline{0000000000110011}$

$\{\frac{1-(-1)^{C(j,11)}}{2}\}_{j=0}^{\infty}=\overline{0000000000010001}$

$\{\frac{1-(-1)^{C(j,12)}}{2}\}_{j=0}^{\infty}=\overline{0000000000001111}$

$\{\frac{1-(-1)^{C(j,13)}}{2}\}_{j=0}^{\infty}=\overline{0000000000000101}$

$\{\frac{1-(-1)^{C(j,14)}}{2}\}_{j=0}^{\infty}=\overline{0000000000000011}$
\end{lemma}
\begin{theorem}\label{th4}
If $t,l$ are positive integers, then
$X(2^t,2^{t+1}l-1)$ is balanced.
\end{theorem}
\begin{proof}
First, $C(j,2^t)=0$ when $0\le j\le 2^t-1$. By Lucas' Theorem, we have
\[
 C(j,2^t)\equiv 1 \pmod{2}\ \text{when}\ 2^t\le j\le 2^{t+1}-1.
 \]
 By Lemma \ref{le8}, the period of $\{(-1)^{C(j,2^t)}\}_{j=0}^{\infty}$ is $2^{t+1}$.
 Hence, we get the sequence $\{(-1)^{C(j,2^t)}\}_{j=0}^{2^{t+1}l-1}$ by repeating
 $\underbrace{++\cdots+}_{2^t}\underbrace{--\cdots-}_{2^t}$ exactly $l$ times.
 Obviously $\{(-1)^{C(j,2^t)}\}_{j=0}^{2^{t+1}l-1}$ is a (trivial) solution of the
 equation $\sum_{i=0}^{n}x_iC(n,i)=0$ when $n=2^{t+1}l-1$. Using Theorem \ref{th2}
 we obtain our result.
\end{proof}

We conjecture that the functions in Theorem \ref{th4} are the only balanced ones.
{\bf Conjecture 1.}
{\em
 There are no nonlinear balanced elementary symmetric polynomials except
for $X(2^t,2^{t+1}\ell-1)$, where $t$ and $\ell$ are any positive integers.
}

\section{Results Concerning Conjecture 1}

The remainder of the paper will be devoted to the study of Conjecture 1.
A Boolean function $f(\xx)$ in $n$ variables is said to satisfy the {\em Strict Avalanche
Criterion} (``is SAC'' for short) if changing any one of the $n$ bits in the input
$\xx$ results in the
output of the function being changed for exactly half of the $2^{n-1}$ vectors $\xx$ with
the changed input bit.
The SAC concept is relevant for our work because of
\begin{lemma}
\label{lem9}
The function  $f(\xx) = X(d,n)$ is SAC if and only if $X(d-1,n-1)$ is balanced.
\end{lemma}
\begin{proof}
  By definition, $f$ is SAC if and only if
\[
               f(\xx) \+ f(\xx\+{\bf a})\ \text{is balanced for all}\ {\bf a} \in  GF(2)^n,\
               \text{with $wt({\bf a}) = 1$}.
\]
 We have $f(\xx) \+ f(\xx\+ (0,\ldots, 0, 1)) = X(d-1,n-1)$,
so the lemma is proved.
\end{proof}

We previously mentioned that
 any symmetric function is completely determined by the weight of its input, that is,
$f(\xx)=v_f(wt(\xx))$. Moreover,
recall the usual algebraic normal form (ANF) of a Boolean function  $f$ in $n$ variables
\[
f(x_1,\ldots,x_n)=\bigoplus_{i=0}^n \lambda_f(i)\bigoplus_{\uu,wt(\uu)=i}
\prod_{j=1}^n x_j^{u_j},
\]
where $\displaystyle v_f(i)=\bigoplus_{j\preceq i} \lambda_f(j)$, and
$\displaystyle \lambda_f(i)=\bigoplus_{j\preceq i} v_f(j)$, over $GF(2)$
($j\preceq i$ means that the binary
expansion of $j$ is less than the binary expansion of $i$, in lexicographical order)
(see \cite[Propositions 1 and 2, p. 2792]{ann}).

The ANF of a symmetric function becomes
\begin{equation}
\label{ANF-sym}
f(x_1,\ldots,x_n)=\bigoplus_{d=0}^n \lambda_f(d) X(d,n),
\end{equation}
in our notations.
Further, when $f$ is an elementary symmetric function, then $\lambda_f(d)=1$ is the only nonzero
coefficient in the representation (\ref{ANF-sym}). Moreover,
\begin{equation}
\label{eqold8}
v_f(i)=\bigoplus_{j\preceq i} \lambda_f(j)=
\begin{cases}
\lambda_f(d),\ if\ d\preceq i\\
0,\ \ otherwise.
\end{cases}
\end{equation}

We need the following further lemmas. We define the well known Walsh transform $W_f(\ww)$ by
\[
 W_f(\ww) = \sum_{\xx \in \BBZ_2^n} (-1)^{f(\xx) +  \xx \cdot \ww}.
\]
\begin{lemma}
\label{lemnew10}
A Boolean function $f$ in $n$ variables is SAC if and only if for every vector
$\bf u$ with $wt(\uu) = 1$
and every vector $\bf v$, we have
\[
        \sum_{\ww \preceq \bar {\uu}}  W_f(\ww \+ \vv)^2 = 2^{wt(  \bar{\uu}) + n}.
\]
\end{lemma}
\begin{proof}
 This is a special case of Proposition 1 of Carlet \cite[p. 35]{Ca99}.
 \end{proof}

\begin{lemma}
\label{lemnew11}  If $f(\xx)$ in $n$ variables is SAC, then
\begin{equation}
\label{eq:new2}
   \sum_{\ww: w_n = 0}  W_f(\ww)^2  =  \sum_{\ww: w_n = 1}  W_f(\ww)^2 = 2^{2n-1}.
\end{equation}
\end{lemma}
\begin{proof}
  We use Lemma \ref{lemnew10} with
$\vv = \bf 0$  and $\uu = (0,\ldots, 0, 1)$.  It follows that $wt(\bar {\uu}) = n -1$,
so the first sum
in (\ref{eq:new2}) equals $2^{2n-1}$.  The two sums add up to $2^{2n}$ by Parseval's
Theorem, so the second sum is also $2^{2n-1}$.
\end{proof}
\begin{lemma}
\label{oldeq4-new12}
 If $f(\xx) = X(d, n)$ is SAC and $d$  is odd, then
\begin{equation}
\label{eqold4}
 W_f({\bf 0}) = 2^n - 2\,wt(f)\  \text{ and}\    W_f({\bf 1}) = 2\,wt(f).
\end{equation}
\end{lemma}
\begin{proof}
 The first equation in (\ref{eqold4}) is clear for any $f$, whether or not $d$ is odd.

 For the second equation, we observe that
by (\ref{eqold8}) our hypotheses imply that $v_f(k) = 0$ for all even $k$.  Since
\[
W_f({\bf 0}) = \sum_{k=0}^n  (-1)^{v_f(k)} C(n,k)\  \text{and}\
W_f({\bf 1}) = \sum_{k=0}^n  (-1)^{v_f(k) + k} C(n,k),
\]
a computation gives
\[
W_f({\bf 0})  +   W_f({\bf 1}) = 2^n.
\]
Now the second equation in (\ref{eqold4}) follows from the first one.
\end{proof}

We define
\begin{equation}
\begin{split}
A&=0,0,1,1;\ {\bar A}=1,1,0,0;\ B=0,1,0,1;{\bar B}=1,0,1,0; \\
C&=0,1,1,0;\ {\bar C}=1,0,0,1;\ D=0,0,0,0;{\bar D}=1,1,1,1.
\end{split}
\end{equation}
     The next two lemmas are used in the proof of our Theorem~\ref{Theorem 1}.
\begin{lemma} {\em (Folklore Lemma [22, Lemma 3.7.2])}
\label{lemfolkbased}
Any affine function
$f$ on $n$ variables, $n\ge 2$, is a linear string of length $2^n$ made up of $4$-bit
blocks $I_1,\ldots,I_{2^{n-2}}$ given as follows:
\begin{enumerate}
\item[1.]   The first block $I_1$ is one of $A,B,C,D,\bar A,\bar
B,\bar C\ {\text{or}}\ \bar D$.
\item[2.]  The second block $I_2$ is $I_1$ or $\bar I_1$.
\item[3.]   The next two blocks $I_3$, $I_4$ are $I_1$, $I_2$ or $
\bar I_1$, $\bar I_2$.
\item[]   \qquad $\cdots \cdots \cdots \cdots \cdots \cdots \cdots
\cdots \cdots \cdots \cdots \cdots $
\item[$n-1$.]   The $2^{n-3}$ blocks $I_{2^{n-3}+1},\ldots
,I_{2^{n-2}}$ are $I_1,\ldots ,I_{2^{n-3}}$ or $\bar I_1,...,\bar I_{2^{n-3}}$.
\end{enumerate}
\end{lemma}
\begin{lemma}
\label{Krawtchouk}
We have  $\displaystyle
\sum_{\xx, wt(\xx)\ even} (-1)^{\xx\cdot \ww} = 0$ for all $\ww \neq {\bf 0}$ or $\bf 1$.
\end{lemma}
\begin{proof}
Let $E(\ww)$ denote the $2^{n-1}$-vector of bits $\xx \cdot \ww \pmod 2$, where
$\bf x$ runs through the $n$-vectors $\bf x$ of even weight in lexicographical order.
Thus $E(\ww)$ lists the exponents in the sum in the lemma.  Consider the $2^{n-1}$
by $n$ array of the vectors $\bf x$ with even weight, taken in lexicographical order.
By the Folklore Lemma, each column in this array is a $2^{n-1}$-vector which gives the truth
table of a nonconstant linear function in $n-1$ variables.  In fact, taking the columns
left to right, the functions are simply $x_1, x_2,\ldots, x_{n-1}$,
$x_1 \+ x_2 \+ \cdots \+ x_{n-1}$.
The vector sum of any subset of at least one and at most $n-1$ of the $n$ columns
(corresponding to $\ww \neq \bf 0$ or $\bf 1$) is thus the truth table of a nonconstant
linear function and so it is balanced.  Each vector $E(\ww)$ is one of these vector sums,
so the sum in the lemma is 0.
\end{proof}

\begin{remark}
The sum in Lemma~\ref{lemfolkbased} is the sum of the Krawtchouk polynomials
\cite[pp. 130 and 150--153]{MacW-Sloane} (variable $y = wt(\ww)$)
\[
  P_k(y,n) = \sum_{\bf x, wt(\xx) =k} (-1)^{\xx \cdot \ww} =
          \sum_{j=0}^k (-1)^j C(y,j) C(n-y, k-j)
\]
of even degree $k$ in $y$.
\end{remark}

\begin{theorem}
\label{Theorem 1}
  If $f(\xx) = X(d,n)$ has odd degree $d$, then $W_f(\ww) = -W_f(\bar {\ww})$
 for all $\ww \neq {\bf 0}$ or $\bf 1$.
 \end{theorem}

\begin{proof}
  Let $f$ be an elementary symmetric function of degree $k$, that is $f=X(d,n)$. We
compute the Walsh transform
\begin{equation}
\label{walsh-bar}
\begin{split}
W_f(\overline\ww)
=& \sum_{\xx\in \BBZ_2^n} (-1)^{f(\xx)+\xx\cdot \overline\ww}\\
=& \sum_{\xx\in \BBZ_2^n} (-1)^{f(\xx)+\xx\cdot ({\bf 1}+\ww)}\\
=& \sum_{\xx\in \BBZ_2^n} (-1)^{f(\xx)+wt(\xx)+\xx\cdot \ww}\\
=& \sum_{k=0}^n\ \sum_{\xx,wt(\xx)=k} (-1)^{f(\xx)+wt(\xx)+\xx\cdot \ww}\\
=& \sum_{k=0}^n(-1)^{v_f(k)+k} \sum_{\xx,wt(\xx)=k} (-1)^{\xx\cdot \ww}.
\end{split}
\end{equation}

Next, we use (\ref{eqold8}). Since $d$ is odd, then any integer $i$ with
$d\preceq i$ has to be odd, as well.
It follows that $v_f(k)=0$, for any even integer $k$.
Thus, (\ref{walsh-bar}) becomes
\begin{eqnarray*}
W_f(\overline\ww)&=&
\sum_{k=0}^n(-1)^{v_f(k)+k} \sum_{\xx,wt(\xx)=k} (-1)^{\xx\cdot \ww}\\
&=&
\sum_{k=0,\,even}^n(-1)^{v_f(k)} \sum_{\xx,wt(\xx)=k} (-1)^{\xx\cdot \ww}\\
&&-
\sum_{k=0,\,odd}^n(-1)^{v_f(k)} \sum_{\xx,wt(\xx)=k} (-1)^{\xx\cdot \ww}\\
&=&  \sum_{\xx,\,wt(\xx)=even} (-1)^{\xx\cdot \ww}-
\sum_{k=0,\,odd}^n(-1)^{v_f(k)} \sum_{\xx,\,wt(\xx)=k} (-1)^{\xx\cdot \ww}.
\end{eqnarray*}
Since
\begin{eqnarray*}
W_f(\ww)
&=&\sum_{k=0,\,even}^n(-1)^{v_f(k)} \sum_{\xx,wt(\xx)=k} (-1)^{\xx\cdot \ww}\\
&& +\sum_{k=0,\,odd}^n(-1)^{v_f(k)} \sum_{\xx,wt(\xx)=k} (-1)^{\xx\cdot \ww}\\
&=&  \sum_{\xx,\,wt(\xx)=even} (-1)^{\xx\cdot \ww}+
\sum_{k=0,\,odd}^n(-1)^{v_f(k)} \sum_{\xx,\,wt(\xx)=k} (-1)^{\xx\cdot \ww},
\end{eqnarray*}
to prove Theorem~\ref{Theorem 1} it will suffice to show that
\[
\sum_{\xx,wt(\xx)=even} (-1)^{\xx\cdot \ww}=0,
\]
as long as $\ww\not= {\bf 0},{\bf 1}$,
and that follows from Lemma~\ref{Krawtchouk}.
\end{proof}

 \begin{theorem}
\label{Theorem 2}
  If $f(x) = X(d,n)$ is SAC and $d$ is odd, then $W_f({\bf 0}) = W_f({\bf 1})$.
  \end{theorem}
  \begin{proof}
By Theorem \ref{Theorem 1}, all of the terms except $W_f(\bf 0)^2$  and
$W_f(\bf 1)^2$
in the two sums in (\ref{eq:new2})
cancel out (for all other $\ww$, $W_f(\ww)$ is in one sum and
$W_f(\bar {\ww})$ is in the other sum).  By Lemma \ref{oldeq4-new12},
both square roots are positive and we get Theorem \ref{Theorem 2}.
  \end{proof}
\begin{corollary}
\label{newcor}
 If $d$ is odd and $f(\xx) = X(d,n)$ is SAC, then $wt(f) = 2^{n-2}$.
 \end{corollary}

Now we determine when $X(d,n)$ is SAC.
To deal with the case when $d$ is an even integer, by Lemma \ref{lem9}, it is enough to show:
\begin{lemma}
\label{case-even}
  If $d > 1$ is odd, then $X(d,n)$ is not balanced.
\end{lemma}
\begin{proof}
Formula (\ref{eqold8}) shows that when $f = X(d,n)$ we have $v_f(i) = 1$ if and
only if $d \preceq i$.  Thus we have
\begin{equation}
\label{eq:new8}
  wt(X(d,n)) = \sum_{d \preceq i, i \leq n}  C(n, i) \leq
\sum_{i\, odd} C(n,i) = 2^{n-1},
\end{equation}
where the inequality holds because $d \preceq i$ and $d$ odd implies $i$ is
odd.  If
$d > 1$, then $d \preceq i$ cannot hold for all odd $i \leq n$
(in particular, $d \not\preceq d - 2$), so the inequality in (\ref{eq:new8}) is strict.
Therefore,
$X(d,n)$ is not balanced.
\end{proof}
\begin{lemma}
\label{lem13}
Suppose $d > 1$ is odd.  If
\begin{equation}
\label{eq9}
d = 2^t + 1\ \text{and}\ n = 2^{t+1} \ell\ \text{for some positive integers}\ t, \ell,
\end{equation}
then $wt(X(d,n)) = 2^{n-2}$.
\end{lemma}
\begin{proof}
 First we observe
\begin{equation}
\label{eq10}
      wt(X(d,n)) = \sum_{d \preceq i, i \leq n}  C(n, i)
\end{equation}
because of (\ref{eqold8}),
which shows that when $f = X(d,n)$ we have $v_f(i) = 1$ if and only if
$d \preceq i$.  By (\ref{eq10}), we need to show that
\begin{equation}
\label{eq11}
 wt(X(d,n)) = \sum_{d \preceq i,\, i \leq n}  C(n, i) = 2^{n-2}
\end{equation}
if and only if (\ref{eq9}) holds.
     If (\ref{eq9}) holds, the sum in (\ref{eq11}) is
\begin{eqnarray*}
&& \sum_{ {2^t+1 \preceq i,\,  i \leq 2^{t+1} \ell} }  C(2^{t+1} \ell, i)  =\\
&& \sum_{2^t+1 \preceq i,\, i \leq 2^{t+1} \ell }
                           (C(2^{t+1} \ell - 1, i) +  C(2^{t+1} \ell - 1, i - 1)) =\\
&& \sum_{2^t \preceq i - 1,\, i - 1 \leq 2^{t+1} \ell - 1}
                           (C(2^{t+1} \ell - 1, i) +  C(2^{t+1} \ell - 1, i - 1)) =\\
&& \sum_{2^t \preceq j,\, j \leq 2^{t+1} \ell - 1}  C(2^{t+1} \ell - 1, j) = 2^{n-2},
\end{eqnarray*}
(note $i$ is never even in the first three sums, since then
$2^t+1 \preceq i$ is false; this justifies the second last equality,
since in the last sum $j$ runs through disjoint pairs of consecutive
integers) where the last sum is $wt(X(2^t, 2^{t+1} \ell - 1)$ by (\ref{eq10})
and so is $2^{n-2}$ by Theorem~\ref{th4}.  Thus we have proved that (\ref{eq9})
implies (\ref{eq11}).
\end{proof}

We would like to prove the converse of the previous lemma.  The following work moves toward
that goal, but does not achieve it.  Next, we prove five lemmas,
which establish many cases of the converse of Lemma~\ref{lem13}.
\begin{lemma}
\label{lem14-case2}
Let $n = 2^{t+1} \ell$ for some positive integers $t,\ell$.
If $j$ is odd and $2^t+1<j<2^{t+1}+1$, then $wt(X(j,n))< 2^{n-2}$.
\end{lemma}
\begin{proof}
     The argument of the previous lemma shows that if (\ref{eq9}) and (\ref{eq11})
     hold for some
given $t$ and $\ell$, then the set
\[
      S(t, \ell) = \{i:\, 2^t + 1 \preceq i,\, i \leq 2^{t+1} \ell = n\}
\]
gives a set of binomial coefficients $\{C(n, i):\, i \in S(t, \ell)\}$ whose
sum is $2^{n-2}$.  (It is easy to see that $S(t, \ell)$ has $n/4$ elements,
but we do not need this fact.)  Now suppose that (\ref{eq11}) holds for
$n = 2^{t+1} \ell$ and for some odd $d = j$, say, satisfying
$2^t + 1 < j < 2^{t+1}+1$.  Then $wt(j) > 2$, so the set
\[
       T(j, n) = \{i:\, j \preceq i,\, i \leq  2^{t+1} \ell = n\}
\]
is a proper subset of  $S(t, \ell)$.  Therefore the sum of the binomial
coefficients in $\{C(n,i):\, i \in T(j, n)\}$ is $< 2^{n-2}$, contradicting
our assumption that (\ref{eq11}) holds with $d = j$.
%
\end{proof}

Since we refer to it often, we include here for completeness an equation given by Canteaut and
Videau in \cite{ann} (these sums are called
{\em lacunary sums of binomial coefficients}, see \cite{Len}). Results like this concerning the binomial
     coefficients are very old.  Some proofs and references are given in \cite{HA76}.
\begin{lemma}
\label{lacunary-lem}
For positive integers $i,n,p$, we have
\begin{equation}
\label{lacunary}
A_n^{2^p}(i)=\sum_{\substack{0\le j\le n\\ j\equiv i\pmod {2^p}}} C(n,j)=
2^{n-p}+2^{1-p}\sum_{j=1}^{2^{p-1}-1} \left(2\cos\left(\frac{j\pi}{2^p} \right)  \right)^n
\cos\left(\frac{j(n-2i)\pi}{2^p}\right)
\end{equation}
\end{lemma}

\begin{lemma}
\label{xia-lem}
Let $t,r$ be positive integers. Suppose that $a_1>a_3\ge a_5\ge \cdots \ge a_J$, with $J=2K+1$,
are nonnegative integers. Define the sum
\[
{\bf T}=\sum_{1\le j\le J} a_j \sin \left(\frac{jr\pi}{2^{t+1}} \right).
\]
If ${\bf T}=0$, then $r\equiv 0\pmod {2^{t+1}}$.
\end{lemma}
\begin{proof}
Write $b_k=a_j$, for $j=2k+1$. For convenience, let $\alpha=\frac{r\pi}{2^{t+1}}$. Then, using Abel's summation formula,
$\bf T$ becomes
\begin{eqnarray*}
{\bf T}
&=&\sum_{k=0}^K b_k\sin ((2k+1)\alpha)\\
&=& \sum_{m=0}^{K-1} (b_m-b_{m+1}) \sum_{k=0}^m \sin ((2k+1)\alpha) + b_K\sum_{k=0}^K \sin ((2k+1)\alpha).
\end{eqnarray*}
Note that for the first term where $m=0$, we have $(b_0-b_1)\sin \alpha\neq 0$, if $r\neq 0\pmod {2^{t+1}}$. Also,
$(b_m-b_{m+1})\ge 0$, and $b_K\ge 0$. The conclusion follows once we show that
\[
\sin\alpha \quad\text{and}\quad \sum_{k=0}^m \sin((2k+1)\alpha)
\]
have the same sign. Indeed
\begin{eqnarray*}
\sin\alpha \sum_{k=0}^m \sin ((2k+1)\alpha)
&=&\frac{1}{2} \sum_{k=0}^m (\cos (2k\alpha)-\cos((2k+2)\alpha)\\
&=& \frac{1}{2} (1-\cos((2m+2)\alpha)\ge 0.
\end{eqnarray*}
The lemma is proved.
\end{proof}

\begin{remark}
\label{rem-xia}
Note that $\bf T$ above has the same sign as $\sin\alpha$.
\end{remark}

Because of Theorem~\ref{th3},
     there is no loss of generality in taking  $n \geq 2(d - 1)$ in our next lemma.
\def\ii{{\bf i}}
\begin{lemma}
\label{lem14-case4}
Let $r,t$ be positive integers, $d = 2^t + 1$, $n = 2^{t+1}$, and $r\not\equiv 0 \pmod {2^{t+1}}$.
Then $wt(X(d,n+r)) \neq 2^{n+r-2}$.
\end{lemma}
\begin{proof}
 Let $d:=1+2^t$ be fixed.
Now, using Pascal's identity, we get that $S:=wt(X(d,n+r))$ satisfies
\begin{eqnarray}
S&=&\sum_{d\preceq i\le n+r} C({n+r},{i})
= \sum_{d\preceq i\le n+r} \left(C({n+r-1},{i})+C({n+r-1},{i-1})\right)\nonumber\\
&=& \sum_{d\preceq i\le n+r-1} C({n+r-1},{i}) +
\sum_{\substack {2^t \preceq j\le n+r-1\\ j\ even}} C({n+r-1},{j})\nonumber\\
&=& \sum_{d\preceq i\le n+r-2} C({n+r-2},{i}) +
\sum_{\substack {2^t \preceq j\le n+r-1\\ j\ even}}
\left(C({n+r-1},{j})+C({n+r-2},{j})\right)\nonumber
\end{eqnarray}
Continuing in this manner, we obtain
\begin{eqnarray}
S&=& \sum_{d\preceq i\le n+r-r} C({2^{t+1}},{i}) +
\sum_{\substack {2^t \preceq j\le n+r-1\\ j\ even}}
\sum_{k=1}^r C({n+r-k},{j})\nonumber\\
&=& 2^{n-2}+ \sum_{\substack {2^t \preceq j\le n+r-1\\ j\ even}}
\sum_{k=1}^r C({n+r-k},{j})\nonumber\\
&=& 2^{n -2}+ \sum_{k=1}^r \sum_{\substack {2^t \preceq j\le n+r-1\\ j\ even}}
C({n+r-k},{j})\nonumber\\
&=& 2^{n -2}+ \sum_{k=1}^r \sum_{s=0}^{2^{t-1}-1}
\sum_{\substack{j\equiv 2s+2^t\pmod {2^{t+1}}\\
0\le j\le n+r-1}} C({n+r-k},{j})
\end{eqnarray}
We push further the previous identity, by computing
the innermost sum. So,
\begin{eqnarray*}
\sum_{\substack{j\equiv 2s+2^t\pmod {2^{t+1}}\\ 0\le j\le n+r-1}}
C({n+r-k},{j})=A_N^{2^{t+1}} (2s+2^t)
\end{eqnarray*}
in the notations of Lemma \ref{lacunary-lem}, where $N:=n+r-k$.
Thus, using equation (\ref{lacunary}), we obtain
\begin{eqnarray*}
A_N^{2^{t+1}} (2s+2^t)=2^{n+r-k-t-1}+2^{-t} \sum_{a=1}^{2^t-1}
\left(2\cos \frac{a\pi}{2^{t+1}}\right)^N
\cos \frac{a(N-4s-2^{t+1})\pi}{2^{t+1}}.
\end{eqnarray*}
Since
\[
\cos \frac{a(N-4s-2^{t+1})\pi}{2^{t+1}}=(-1)^a \cos \frac{a(N-4s)\pi}{2^{t+1}},
\]
we get
\begin{equation}
A_N^{2^{t+1}} (2s+2^t)=2^{n+r-k-t-1}+2^{-t} \sum_{a=1}^{2^t-1} (-1)^a
\left(2\cos \frac{a\pi}{2^{t+1}}\right)^N
 \cos \frac{a(N-4s)\pi}{2^{t+1}}
\end{equation}
 We obtain
\begin{eqnarray*}
S&=&2^{n-2}+ \sum_{k=1}^r \sum_{s=0}^{2^{t-1}-1} A_N^{2^{t+1}} (2s+2^t)\\
&=& 2^{n-2}+\sum_{k=1}^r \sum_{s=0}^{2^{t-1}-1}2^{n+r-k-t-1}\\
&&+2^{-t}\sum_{k=1}^r \sum_{s=0}^{2^{t-1}-1} \sum_{a=1}^{2^t-1}
(-1)^a\left(2\cos \frac{a\pi}{2^{t+1}}\right)^N
 \cos \frac{a(N-4s)\pi}{2^{t+1}}\\
 &=& 2^{n-2}+2^{n+r-2}\sum_{k=1}^r 2^{-k}+
 2^{-t}\sum_{k=1}^r \sum_{s=0}^{2^{t-1}-1} \sum_{a=1}^{2^t-1}
 (-1)^a\left(2\cos \frac{a\pi}{2^{t+1}}\right)^N
 \cos \frac{a(N-4s)\pi}{2^{t+1}}\\
 &=&2^{n+r-2}+
 2^{-t}\sum_{k=1}^r \sum_{s=0}^{2^{t-1}-1} \sum_{a=1}^{2^t-1}
 (-1)^a\left(2\cos \frac{a\pi}{2^{t+1}}\right)^N
 \cos \frac{a(N-4s)\pi}{2^{t+1}}.
\end{eqnarray*}
Therefore, to prove our assertion, we need to show that
\begin{eqnarray*}
T:&=&\sum_{k=1}^r \sum_{s=0}^{2^{t-1}-1} \sum_{a=1}^{2^t-1} (-1)^a
\left(2\cos \frac{a\pi}{2^{t+1}}\right)^{n+r-k}
 \cos \frac{a(n+r-k-4s)\pi}{2^{t+1}}\\
&=&\sum_{k=1}^r \sum_{a=1}^{2^t-1} (-1)^a\left(2\cos \frac{a\pi}{2^{t+1}}\right)^{n+r-k}
\sum_{s=0}^{2^{t-1}-1}\cos \frac{a(n+r-k-4s)\pi}{2^{t+1}}\neq 0.
\end{eqnarray*}
Since
\[
\frac{a(n+r-k-4s)\pi}{2^{t+1}}=a\pi+\frac{(r-k-4s)a\pi}{2^{t+1}},
\]
and so,
\[
\cos\left(\frac{a(n+r-k-4s)\pi}{2^{t+1}}\right)=(-1)^{a}
\cos\left(\frac{(r-k-4s)a\pi}{2^{t+1}}\right),
\]
we obtain
\[
T=\sum_{k=1}^r \sum_{a=1}^{2^t-1}
\left(2\cos \frac{a\pi}{2^{t+1}}\right)^{n+r-k}
\sum_{s=0}^{2^{t-1}-1}  \cos\left(\frac{(r-k-4s)a\pi}{2^{t+1}}\right)
\]

Formula (17.1.1) of \cite{Hansen} states
\begin{equation}
\label{eq:star1}
\sum_{s=0}^N \cos(sx+y)=\csc \frac{x}{2} \cos\left(\frac{Nx}{2}+y\right)
\sin\left(\frac{(N+1)x}{2}\right).
\end{equation}
Taking $A=\frac{a\pi}{2^{t+1}}$, $N=2^{t-1}-1$, $x=-4A$, $y=(r-k)A$ in the previous formula,
we obtain
\begin{eqnarray*}
\sum_{s=0}^{2^{t-1}-1} \cos\left((r-k-4s)A\right)
&=& \csc (-2A) \cos\left((2^{t-1}-1) (-2A)+(r-k)A\right)\sin (2^{t-1}(-2A))\\
&=& \csc(2A) \sin\left(\frac{a\pi}{2}\right)\cos\left(-\frac{a\pi}{2}+(r-k+2)A\right)\\
&=& \frac{1-(-1)^a}{2} \frac{\sin((r-k+2)A)}{\sin(2A)}.
\end{eqnarray*}
Now, $T$ becomes
\begin{eqnarray*}
T&=&\sum_{k=1}^r \sum_{a=1}^{2^t-1}  \frac{1-(-1)^a}{2}
\left(2\cos A\right)^{n+r-k}  \frac{\sin((r-k+2)A)}{\sin(2A)}\\
&=& \sum_{a=1}^{2^t-1}  \frac{1-(-1)^a}{2}
\frac{\left(2\cos A\right)^{n+r}}{\sin(2A)}
\sum_{k=1}^r \left(2\cos A\right)^{-k}  \sin((r-k+2)A)
\end{eqnarray*}
We evaluate the inside sum using formula (14.7.1) of \cite{Hansen}
\begin{eqnarray*}
\sum_{k=1}^{N-1} b^k \sin(kx+y)&=&-\sin y+(1-2b\cos x+b^2)^{-1}\cdot\\
&&[\sin y+b\sin(x-y)-b^N \sin(Nx+y)\\
&&+b^{N+1}\sin ((N-1)x+y)]
\end{eqnarray*}
with $N=r+1$, $b=(2\cos A)^{-1}$, $x=-A$, $y=(r+2)A$. We get
\begin{eqnarray*}
&&\sum_{k=1}^r (2\cos A)^{-k} {\sin((r-k+2)A)}\\
&&= -\sin((r+2)A) +b^{-2}(\sin ((r+2)A)-b  \sin ((r+3)A) \\
&&\qquad \qquad -b^{r+1}\sin A+b^{r+2}\sin(2A))\\
&&= -\sin((r+2)A) +b^{-1}(2\cos A\sin((r+2)A)-\sin((r+3)A))\\
&&\qquad\qquad  -b^r(2\cos A\sin A-\sin(2A))\\
&&= -\sin((r+2)A) +2\cos A\sin ((r+1)A)=\sin (rA).
\end{eqnarray*}
and so,
\begin{equation*}
\begin{split}
T&=\sum_{a=1}^{2^t-1}  \frac{1-(-1)^a}{2}
\left(2\cos A\right)^{n+r-1}\frac{\sin (rA)}{\sin A}\\
&=  \sum_{a=1,\ odd}^{2^t-1}
\left(2\cos A\right)^{n+r-1}
\frac{\sin (rA)}{\sin A}
\end{split}
\end{equation*}

Recall that our initial sum is
\[
S=2^{n+r-2}+2^{-t} T,
\]
so we need to prove $T \neq 0$.
Observing that
\[
a_j=\left(\cos\frac{j\pi}{2^{t+1}} \right)^{2^{t+1}+r-1}\cdot \frac{1}{\sin\frac{j\pi}{2^{t+1}}}
\]
strictly decreases as $j$ increases, $1\le j\le 2^t-1$, Lemma \ref{xia-lem}
shows that $T\neq 0$, thereby proving our claim. (One can prove, by a slightly more complicated method that,
in fact, $T>0$, but we did not need that.)
The proof of the lemma is done.
\end{proof}

\begin{lemma}
\label{lem14-case3}
If $d$ is odd and $2^t + 1<d\le 2^{t+1}-1$ for some positive integer $t$,
then $wt(X(d,n)) \neq 2^{n-2}$ for any $n$ of the form $n = 2^{t+1} \ell+r$,
where $\ell$ is even and $0 \leq r < 2^{t+1} + 2^t$.
\end{lemma}
\begin{proof}
From  equation (\ref{eq10})
we have
\begin{equation}
\label{eqA}
  wt(X(2^t + 1, n)) = \sum_{k \in I(t)}\ \sum_{i \equiv k\, ({\rm mod}\ {2^{t+1}}),\, i \leq n} C(n, i)
\end{equation}
where
\[
\begin{split}
   I(t) &= \{k:\ k\ odd,\ 2^t + 1 \leq k \le 2^{t+1} - 1\} \\
   &=
   \{\text{the largest $2^{t-1}$ odd least positive residues $\pmod{2^{t+1}}$}\}.
\end{split}
\]

Let $k:=2^t+2s+1$, where $0\le s\le 2^{t-1}-1$, and
let $A_n^{2^{t+1}}(k)$ denote the inner sum in (\ref{eqA}).
Then Lemma \ref{lacunary-lem} gives (with $A=\frac{j \pi}{2^{t+1}}$)
\begin{equation}
\label{eqB}
\begin{split}
A_n^{2^{t+1}}(k)
&= 2^{n-(t+1)} + 2^{n-t}\sum_{j=1}^{2^t-1}
 (\cos A)^n
 \cos ((n-2k) A)\\
&= 2^{n-(t+1)} + 2^{n-t}\sum_{j=1}^{2^t-1} (-1)^j
 (\cos A)^n
\cos\left((n-2-4s) A\right),
 \end{split}
\end{equation}
since
\begin{eqnarray*}
\cos ((n-2k) A)&=&\cos ((n-2(2^t+2s+1)) A)\\
&=&\cos ((n-4s-2) A-2^{t+1} A)= \cos ((n-4s-2) A- j\pi)\\
&=&\cos((n-4s-2)A)\cos(j\pi)+\sin((n-4s-2)A) \sin(j\pi)\\
&=&(-1)^j\cos((n-4s-2)A).
\end{eqnarray*}

If $d$ is odd, 
let $J(d)\subset I(t)$ be the subset of $I(t)$, made up of the $2^{t-2}$ integers $k$
that satisfy $d\preceq k\le 2^{t+1}-1$
 (for example, if $d = 2^t + 3$,
then $J(d)$ contains every other integer in $I(t)$, starting with $2^t + 3$).
Let $n=2^{t+1}\ell+r$, $0\le r<2^{t+1}+2^t$. If $r=0$, Lemma \ref{lem14-case2} implies the result.
Now, assume $1\le r< 2^{t+1}+2^t.$
Using (\ref{eq:star1}) we obtain (recall that $A=\frac{j \pi}{2^{t+1}}$)
\begin{equation}
\label{eq:sum-cos2}
\begin{split}
\sum_{s=0}^{2^{t-1}-1} \cos\left(s(-4A)+(n-2)A \right)
&=\csc(-2A)\cos\left((2^{t-1}-1)(-2A)+(n-2)A\right)\sin(2^{t-1}(-2A))\\
&= \csc (2A) \cos(-2^t A+nA) \sin\left(\frac{j\pi}{2}\right)\\
&=\csc(2A) \left(\cos\left(\frac{j\pi}{2}\right)\cos (nA)+\sin\left(\frac{j\pi}{2}\right)\sin (nA)
 \right)\sin\left(\frac{j\pi}{2}\right)\\
 &=\csc(2A)\sin^2\left(\frac{j\pi}{2}\right) \sin(nA)\\
&=\frac{1-(-1)^j}{2} \csc(2A) \sin((2^{t+1}\ell+r)A)\\
&= \frac{1-(-1)^j}{2} \csc(2A) (-1)^\ell \sin (rA).
\end{split}
\end{equation}

Certainly (with $k=2^t+2s+1$),
\begin{eqnarray*}
 wt(X(d,n))
&=& \sum_{k \in J(d)}\ \sum_{i \equiv k\, ({\rm mod}\ {2^{t+1}}),\, i \leq n} C(n, i)\\
&\le&
\sum_{k \in I(t)} A_n^{2^{t+1}}(k)=\sum_{s=0}^{2^{t-1}-1} A_n^{2^{t+1}}(2^t+2s+1).
\end{eqnarray*}
Then, using (\ref{eqB}) and (\ref{eq:sum-cos2})
\begin{equation}
\label{eqD}
\begin{split}
 & \sum_{s=0}^{2^{t-1}-1} A_n^{2^{t+1}}(2^t+2s+1)
= 2^{n-2}+2^{n-t}\sum_{s=0}^{2^{t-1}-1}\ \sum_{j=1}^{2^t-1} (-1)^j
 \left(\cos A\right)^n\,
\cos\left((n-2-4s) A\right)\\
&= 2^{n-2}+2^{n-t}\sum_{j=1}^{2^t-1} (-1)^j
 \left(\cos A\right)^n\,
 \sum_{s=0}^{2^{t-1}-1}
\cos\left((n-2-4s)A \right)\\
&= 2^{n-2}+2^{n-t}\sum_{j=1}^{2^t-1} (-1)^{\ell+j}
\left(\cos A\right)^n\,
 \frac{1-(-1)^j}{2} \frac{\sin (rA)}{\sin (2A)}\\
&= 2^{n-2}+2^{-t}(-1)^{\ell+1}\sum_{j=1,odd}^{2^t-1}
\left(2\cos A\right)^{n-1}\,\frac{\sin (rA)}{\sin A}:=S
  \end{split}
\end{equation}
But the last sum is strictly positive by Lemmas~\ref{xia-lem} and \ref{lem14-case4}.
Therefore, if $\ell$ is even, $S<2^{n-2}$, and this proves our lemma.
\end{proof}

\begin{remark}
We see that if $n=2^{t+1}\ell+r$, $\ell$ odd and   $r<2^t$, then we can write
$n=2^{t+1}\ell+r=2^{t+1}(\ell-1) +2^{t+1}+r$, with $\ell-1$ even, and
$0\le r':=2^{t+1}+r < 2^{t+1}+2^t$.
Thus, the only cases left unchecked in the previous lemma (which gives
many cases of Conjecture $1$) are:
 $n = 2^{t+1}\ell + r$, $\ell$ odd, $2^t \leq r < 2^{t+1}$.
\end{remark}

\section{The Case $wt(d)\geq 3$}

 Lemma \ref{lem9}, Corollary \ref{newcor} and Lemma \ref{lem14-case4} show that Conjecture 1
  holds for any $X(d, n)$ with $d = 2^t$.  A key fact, given in the proof of Lemma \ref{lem14-case4},
  is a useful formula for $wt(X(d,n))$ when $wt(d) = 2$.  We can find a similar formula when $wt(d) = 3$,
  however it becomes substantially harder to handle.
\begin{lemma}
\label{lem22}
Let $d:=1+2^{s}+2^t$, where $1\le s < t$ and $t\ge 2$. Then
\begin{equation}
\begin{split}
\label{eq:simpl1}
wt(X(d,n))&=2^{n-3}-2^{-t} \sum_{j=1,odd}^{2^t-1}
(2\cos A)^{n-1} \frac{\sin((n-2^s)A)\sin (2^s A) }{\sin A\sin (2^{s+1}A)}\\
&\quad-2^{-s-1} \sum_{k=1,odd}^{2^s-1} (2\cos B)^{n-1} \frac{\sin(nB)}{\sin B}
\end{split}
\end{equation}
\end{lemma}
\begin{proof}
Let $A=\frac{j\pi}{2^{t+1}}$, $B=\frac{k\pi}{2^{s+1}}$.
From $d\preceq i$, we get that $i=2^{t+1}i'+2^t+2^{s+1}p+2^s+2q+1$, and so, $i\equiv 2^t+2^{s+1}p+2^s+2q+1\pmod {2^{t+1}}$.
Certainly the converse is also true.
Using the previous observation,
\begin{equation}
\begin{split}
&wt(X(d,n))
=\sum_{d\preceq i\le n} C(n,i)=\sum_{p=0}^{2^{t-s-1}-1}\ \sum_{q=0}^{2^{s-1}-1} A_n^{2^{t+1}}(2^t+2^{s+1}p+2^s+2q+1)\\
&=\sum_{p=0}^{2^{t-s-1}-1}\ \sum_{q=0}^{2^{s-1}-1} \left(2^{n-t-1}+2^{-t}\sum_{j=1}^{2^t-1} (2 \cos A)^n
\cos((n-2^{t+1}-2^{s+2}p-2^{s+1}-4q-2)A) \right)\\
&= 2^{t-s-1}2^{s-1} 2^{n-t-1}+2^{-t}\sum_{j=1}^{2^t-1}(2 \cos A)^n\sum_{p=0}^{2^{t-s-1}-1}\ \sum_{q=0}^{2^{s-1}-1}
\cos((n-2^{t+1}-2^{s+2}p-2^{s+1}-4q-2)A)\\
&=2^{n-3}+2^{-t}\sum_{j=1}^{2^t-1}(2 \cos A)^n \sum_{p=0}^{2^{t-s-1}-1}\ \sum_{q=0}^{2^{s-1}-1}
\cos((n-2^{t+1}-2^{s+2}p-2^{s+1}-4q-2)A)
\end{split}
\end{equation}
using Lemma~\ref{lacunary-lem}.
Further, by using formula (\ref{eq:star1}) with
$x=-4A$, $y=(n-2^{t+1}-2^{s+2}p-2^{s+1}-2)A$, $N=2^{s-1}-1$, the innermost sum is equal to
\begin{equation*}
\begin{split}
&\csc(x/2)\cos(Nx/2+y)\sin ((N+1)x/2)\\
&=\csc(-2A)\cos((2^{s-1}-1)(-2A)+(n-2^{t+1}-2^{s+2}p-2^{s+1}-2)A) \sin(2^{s-1}(-2A))\\
&=\csc (2A)\cos((n-2^{t+1}-2^{s+2}p-3\cdot 2^s)A) \sin (2^s A),
\end{split}
\end{equation*}
which is defined everywhere, since $j\le 2^t-1$.
Thus,
\begin{equation}
\label{eq:genX}
\begin{split}
&wt(X(d,n))
=2^{n-3}+2^{-t}\sum_{j=1}^{2^t-1}(2 \cos A)^{n-1} \frac{\sin (2^s A)}{\sin A} \sum_{p=0}^{2^{t-s-1}-1}
\cos((n-2^{t+1}-2^{s+2}p-3\cdot 2^s)A).
\end{split}
\end{equation}
Let
\[
U:=\{j:\, j=2^{t-s}k, 1\le k\le 2^s-1 \}
\]
We distinguish two cases:\\
{\bf Case 1.} Assume $j\in U$. That means that
\[
2^{s+2}A=2^{s+2} \frac{j\pi}{2^{t+1}}=2^{s+2} \frac{k 2^{t-s}\pi}{2^{t+1}}=2k\pi,
\]
and using the periodicity of the cosine function, we obtain that in this case, the innermost sum is
\[
2^{t-s-1} \cos((n-2^{t+1}-3\cdot 2^s)A).
\]

\noindent
{\bf Case 2.} Assume $j\not\in U$. In this case, we apply again formula (\ref{eq:star1}) with
$x=-2^{s+2}A$, $y=(n-2^{t+1}-3\cdot 2^s)A$, $N=2^{t-s-1}-1$, the innermost sum is equal to
\begin{eqnarray*}
&&\csc(-2^{s+1}A)\cos((2^{t-s-1}-1)(-2^{s+1}A)+(n-2^{t+1}-3\cdot 2^s)A)\sin (2^{t-s-1}(-2^{s+1}A))\\
&&=\csc(2^{s+1}A) \cos(-2^t A+(n-2^{t+1}-2^s)A)\sin (2^tA)\\
&&=\csc(2^{s+1}A)\cos((n-2^s)A-3j\pi/2)\sin(j\pi/2)\\
&&=\csc(2^{s+1}A)\cos((n-2^s)A+j\pi/2)\sin(j\pi/2)\
\end{eqnarray*}

Thus, from equation (\ref{eq:genX}), we obtain (note that $A=B$, if $j=2^{t-s}k$; also, $2^{t+1}A=j\pi$, $2^sB=k\pi/2$)
\begin{equation}
\begin{split}
\label{eq:unsimpl}
wt(X(d,n))&= 2^{n-3}+2^{-t} \sum_{j=1,j\not\in U}^{2^t-1}
(2\cos A)^{n-1} \frac{\cos((n-2^s)A+j\pi/2)\sin(j\pi/2)\sin (2^s A) }{\sin A\sin (2^{s+1}A)}\\
&\quad+2^{-t} \sum_{j=1,j\in U}^{2^t-1} (2\cos A)^{n-1} \frac{\sin (2^s A)}{\sin A} 2^{t-s-1} \cos((n-3\cdot 2^s)A-j\pi)\\
&= 2^{n-3}+2^{-t} \sum_{j=1,j\not\in U}^{2^t-1}
(2\cos A)^{n-1} \frac{\cos((n-2^s)A+j\pi/2)\sin(j\pi/2)\sin (2^s A) }{\sin A\sin (2^{s+1}A)}\\
&\quad+2^{-s-1} \sum_{k=1}^{2^s-1} (2\cos B)^{n-1} \frac{\sin(k\pi/2)}{\sin B} \cos((n-3\cdot 2^s)B-2^{t-s}k\pi)\\
&= 2^{n-3}+2^{-t} \sum_{j=1,j\not\in U}^{2^t-1}
(2\cos A)^{n-1} \frac{\cos((n-2^s)A+j\pi/2)\sin(j\pi/2)\sin (2^s A) }{\sin A\sin (2^{s+1}A)}\\
&\quad+2^{-s-1} \sum_{k=1}^{2^s-1} (2\cos B)^{n-1} \frac{\sin(k\pi/2)}{\sin B} \cos(nB+k\pi/2).
\end{split}
\end{equation}
(The last equality follows from the periodicity of $\cos$, and also from
$\cos((n-3\cdot 2^s)B)=\cos(nB-3 k\pi/2)=\cos(nB+k\pi/2)$.)
Further, if $j\not\in U$, then $\sin (2^{s+1}A)$ is well defined, however $\sin (j\pi/2)=0$, if $j$ is even. Thus, the terms in the
first sum of the last equation of (\ref{eq:unsimpl}) are zero, unless $j$ is odd. Then, if $j$ is odd, we get
\begin{eqnarray*}
&&\cos((n-2^s)A+j\pi/2)\sin(j\pi/2)\\
&=&\left(\cos((n-2^s)A) \cos(j\pi/2)-\sin((n-2^s)A)\sin(j\pi/2)\right)\sin(j\pi/2)\\
&=&-\sin((n-2^s)A).
\end{eqnarray*}
Therefore,
\begin{equation*}
\begin{split}
wt(X(d,n))&=2^{n-3}-2^{-t} \sum_{j=1,odd}^{2^t-1}
(2\cos A)^{n-1} \frac{\sin((n-2^s)A)\sin (2^s A) }{\sin A\sin (2^{s+1}A)}\\
&\quad+2^{-s-1} \sum_{k=1}^{2^s-1} (2\cos B)^{n-1} \frac{\sin(k\pi/2)}{\sin B} \cos(nB+k\pi/2)
\end{split}
\end{equation*}
or better, yet,
\begin{equation*}
\begin{split}
wt(X(d,n))&=2^{n-3}-2^{-t} \sum_{j=1,odd}^{2^t-1}
(2\cos A)^{n-1} \frac{\sin((n-2^s)A)\sin (2^s A) }{\sin A\sin (2^{s+1}A)}\\
&\quad-2^{-s-1} \sum_{k=1,odd}^{2^s-1} (2\cos B)^{n-1} \frac{\sin(nB)}{\sin B}
\end{split}
\end{equation*}
\end{proof}

In order to prove Conjecture 1, by Lemma \ref{lem9} and Corollary \ref{newcor} it would suffice to show that
for $n \geq 2(d - 1)$  (we can assume this because of Theorem \ref{th3}) we have
\begin{equation}
\label{eq26}
                                 wt(X(d, n)) \neq  2^{n-2}
\end{equation}
for all pairs $d, n$ except  $d = 2^t + 1, n = 2^{t+1}\ell$, where $t$ and $\ell$ are any positive integers.

     Lemma \ref{lem14-case4}  proves (\ref{eq26}) when $wt(d) = 2$.  We attempted to prove (\ref{eq26})
when $wt(d) = 3$ by using Lemma \ref{lem22}, but the sums in (\ref{eq:simpl1}) were too complicated to
allow us to cover all of the cases.  Certainly (\ref{eq:simpl1}) shows that for fixed $d$, (\ref{eq26})
holds for all sufficiently large $n$, because the factors  $(\cos A)^{n-1}$ and $(\cos B)^{n-1}$ tend to 0 as
$n \to\infty$, which implies  $wt(X(d, n)) - 2^{n-2} < 0$ for all large $n$.  Our computations
suggest that this inequality will always hold if $wt(d)$ is large enough.  In fact, we conjecture

\noindent
{\bf Conjecture 2.}  If $n \geq 2(d - 1)$, $d$ is fixed and $wt(d) \geq 6$, then $wt(X(d, n)) - 2^{n-2} <~0$.

\noindent
{\bf Acknowledgements.}
The authors would like to thank Prof. Jingbo Xia for the proof of Lemma \ref{xia-lem}, which simplified
their original argument.

\end{document}